\documentclass[12pt]{amsart}
\usepackage{amsfonts}
\usepackage{amssymb}

\newcommand{\R}{\mathbb R}

\newcommand{\E}{\mathbb E}

\renewcommand{\span}{\mathrm{span}}
\newcommand{\tr}{\mathrm{tr}}

\newtheorem{thm}{Theorem}[section]
\newtheorem{cor}[thm]{Corollary}
\newtheorem{lem}[thm]{Lemma}
\newtheorem{prop}[thm]{Proposition}
\theoremstyle{definition}

\theoremstyle{remark}

\begin{document}

\title[ON THE THEORY OF SURFACES IN THE FOUR-DIMENSIONAL EUCLIDEAN SPACE]
{ON THE THEORY OF SURFACES IN THE FOUR-DIMENSIONAL EUCLIDEAN SPACE}
\author{Georgi Ganchev and Velichka Milousheva}
\address{Bulgarian Academy of Sciences, Institute of Mathematics and Informatics,
Acad. G. Bonchev Str. bl. 8, 1113 Sofia, Bulgaria}
\email{ganchev@math.bas.bg}
\address{Bulgarian Academy of Sciences, Institute of Mathematics and Informatics,
Acad. G. Bonchev Str. bl. 8, 1113, Sofia, Bulgaria}
\email{vmil@math.bas.bg}

\begin{abstract}

For a two-dimensional surface $M^2$ in the four-dimensional Euclidean space $\E^4$ we
introduce an invariant linear map of Weingarten type in the tangent space of the
surface, which generates two invariants $k$ and $\varkappa$.

The condition $k = \varkappa = 0$ characterizes the surfaces consisting of flat points.
The minimal surfaces are characterized by the equality $\varkappa^2-k=0$.
The class of the surfaces with flat normal connection is characterized by the condition
$\varkappa = 0$. For the surfaces of general type we obtain a geometrically determined
orthonormal frame field at each point and derive Frenet-type derivative formulas.

We apply our theory to the class of the rotational surfaces in $\E^4$,
which prove to be surfaces with flat normal connection, and describe the
rotational surfaces with constant invariants.

\end{abstract}
\maketitle

\section{Introduction} \label{S:Intr}

In \cite{Otsuki1} T. \={O}tsuki introduced
curvatures $\lambda_1, \, \lambda_2, \dots, \lambda_n$ ($\lambda_1
\geq \lambda_2 \geq \dots \geq \lambda_n$) for a surface $M^2$ in
a $(2+n)$-dimensional Euclidean space $\E^{2+n}$,
defining a quadratic form in the normal space of the surface.
In a suitable local frame of the normal space this quadratic form can be
written in a diagonal form and the functions $\lambda_{\alpha}, \, \alpha = 1, \dots, n$
are the coefficients in the diagonalized form
($\lambda_{\alpha}$ is called the \textit{$\alpha$-th curvature} of $M^2$).
These curvatures are closely related to the Gauss curvature $K$ of $M^2$:
$$K = \lambda_1 + \lambda_2 + \dots + \lambda_n.$$
The local cross-section, which diagonalizes the quadratic form is called
a \textit{Frenet cross-section} (\textit{Frenet-frame}) of the surface.

For a surface  $M^2$ in the four-dimensional Euclidean space $\E^4$ the curvatures
$\lambda_1$ and $\lambda_2$ are the maximum and minimum, respectively of the
Lipschitz-Killing curvature of the surface \cite{Otsuki2}. The function $\lambda_1$ is
called the \textit{principal curvature} and the function $\lambda_2$ - the
\textit{secondary curvature} of $M^2$ in $\E^4$.

Using the idea of the Frenet-frames, Shiohama \cite{Shiohama1} proved that a complete
connected orientable surface $M^2$ in  $\E^4$ with curvatures $\lambda_1 = \lambda_2 = 0$
is a cylinder. The same  result is proved in \cite{Shiohama2} for
a surface in a higher dimensional space $\E^{2+n}$.

Our aim is to find invariants of a surface $M^2$ in  $\E^4$,
considering a geometrically determined linear map (of Weingarten type)
in the tangent space of the surface, as well as to obtain a geometric Frenet-type
frame field of $M^2$.

In Section \ref{S:Weingarten} we define a geometrical linear map in the tangent space
of a surface $M^2$ in $\E^4$ and determine a second fundamental form $II$ of the surface.
We find invariants $k$ and $\varkappa$ of $M^2$
(which are analogous to the Gauss curvature and the mean curvature of a surface in $\E^3$).
These invariants divide the points of $M^2$ into four types: flat, elliptic, parabolic and
hyperbolic.

In Section \ref{S:flat} we give a local geometric description of the surfaces
consisting of flat points, proving that they are either planar
surfaces (Proposition \ref{P:3.1}) or developable ruled surfaces (Proposition \ref{P:3.2}).

In Section \ref{S:minimal} we characterize the minimal surfaces in $\E^4$ in terms of the
invariants $k$ and $\varkappa$ (Proposition \ref{P:minimal}).

For the surfaces of general type (which are not minimal and which have no flat points)
in Section \ref{S:general} we obtain a geometrically determined orthonormal frame field
$\{x,y,b,l\}$  at each point of the surface and derive Frenet-type derivative formulas.
The tangent frame field $\{x,y\}$ is determined by the defined second fundamental form $II$,
while the normal frame field $\{b,l\}$ is determined by the mean curvature vector field
of the surface.

We also characterize the surfaces with flat normal connection
in terms of the invariant $\varkappa$ (Theorem \ref{P:flat-normal-1}).

In the last section we apply our theory to the class of the rotational surfaces in $\E^4$,
which prove to be surfaces with flat normal connection, and describe the
rotational surfaces with $k = {\rm const}$.

\section{The Weingarten map} \label{S:Weingarten}

We denote by $g$ the standard metric in the four-dimensional Euclidean space ${\E}^4$
and by $\nabla'$ its flat Levi-Civita
connection. All considerations in the present paper are local and all functions,
curves, surfaces, tensor fields etc. are assumed to be of the class $\mathcal C^{\infty}$.

Let $M^2: z = z(u,v), \, \, (u,v) \in {\mathcal D}$ (${\mathcal D} \subset \R^2$) be a
2-dimensional surface in $\E^4$. The tangent space to $M^2$ at an arbitrary point
$p=z(u,v)$ of $M^2$ is ${\rm span} \{z_u, z_v\}$.

For an arbitrary orthonormal normal frame field $\{e_1, e_2\}$ of $M^2$ we have the
standard derivative formulas:
$$\begin{array}{l}
\vspace{2mm} \nabla'_{z_u}z_u=z_{uu} = \Gamma_{11}^1 \, z_u +
\Gamma_{11}^2 \, z_v
+ c_{11}^1\, e_1 + c_{11}^2\, e_2;\\
\vspace{2mm} \nabla'_{z_u}z_v=z_{uv} = \Gamma_{12}^1 \, z_u +
\Gamma_{12}^2 \, z_v
+ c_{12}^1\, e_1 + c_{12}^2\, e_2;\\
\vspace{2mm} \nabla'_{z_v}z_v=z_{vv} = \Gamma_{22}^1 \, z_u +
\Gamma_{22}^2 \, z_v
+ c_{22}^1\, e_1 + c_{22}^2\, e_2,\\
\end{array} \leqno{(2.1)}$$
where $\Gamma_{ij}^k$ are the Christoffel's symbols and $c_{ij}^k$, $i, j, k = 1,2$
are functions on $M^2$.

We use the standard denotations \;$E(u,v)=g(z_u,z_u), \; F(u,v)=g(z_u,z_v), \;
G(u,v)=g(z_v,z_v)$ for the coefficients of the
first fundamental form and set $W=\sqrt{EG-F^2}$. If $\sigma$
denotes the second fundamental tensor of $M^2$, then we have
$$\begin{array}{l}
\sigma(z_u,z_v)=c_{11}^1\, e_1 + c_{11}^2\, e_2,\\
[2mm]
\sigma(z_u,z_v)=c_{12}^1\, e_1 + c_{12}^2\, e_2,\\
[2mm] \sigma(z_v,z_v)=c_{22}^1\, e_1 + c_{22}^2\,
e_2.\end{array}$$

\vskip 2mm
We introduce the following functions:
$$\Delta_1 = \left|%
\begin{array}{cc}
\vspace{2mm}
  c_{11}^1 & c_{12}^1 \\
  c_{11}^2 & c_{12}^2 \\
\end{array}%
\right|; \quad
\Delta_2 = \left|%
\begin{array}{cc}
\vspace{2mm}
  c_{11}^1 & c_{22}^1 \\
  c_{11}^2 & c_{22}^2 \\
\end{array}%
\right|; \quad
\Delta_3 = \left|%
\begin{array}{cc}
\vspace{2mm}
  c_{12}^1 & c_{22}^1 \\
  c_{12}^2 & c_{22}^2 \\
\end{array}%
\right|;$$
$$L(u,v) = \displaystyle{\frac{2 \Delta_1}{W}, \quad M(u,v) =
\frac{\Delta_2}{W}, \quad N(u,v) = \frac{2 \Delta_3}{W}}.$$

If
$$\begin{array}{l}
\vspace{2mm}
u = u(\bar u,\bar v);\\
\vspace{2mm} v = v(\bar u,\bar v),
\end{array}
\quad (\bar u,\bar v) \in \bar{\mathcal D}, \,\, \bar{\mathcal D} \subset
\R^2 \leqno{(2.2)}$$ is a smooth change of the parameters $\{u,v\}$
on $M^2$ with $J = u_{\bar u} \, v_{\bar v} - u_{\bar v} \,
v_{\bar u}\neq 0$, then
$$\begin{array}{l}
\vspace{2mm}
z_{\bar u} = z_u \,u_{\bar u} + z_v \,v_{\bar u},\\
\vspace{2mm} z_{\bar v} = z_u \,u_{\bar v} + z_v \,v_{\bar v}.
\end{array}$$
Let
$$\begin{array}{l}
\vspace{2mm}
\sigma(z_{\bar u},z_{\bar u}) = \bar c_{11}^1 \, e_1 + \bar c_{11}^2 \, e_2,\\
\vspace{2mm}
\sigma(z_{\bar u},z_{\bar v}) = \bar c_{12}^1 \, e_1 + \bar c_{12}^2 \, e_2,\\
\vspace{2mm} \sigma(z_{\bar v},z_{\bar v}) = \bar c_{22}^1 \, e_1
+ \bar c_{22}^2 \, e_2.
\end{array}$$
Differentiating (2.2) and taking into account (2.1) we find
$$\begin{array}{l}
\vspace{2mm} \bar c_{11}^k = u_{\bar u}^2 \,c_{11}^k  +
2 u_{\bar u}\, v_{\bar u}\, c_{12}^k + v_{\bar u}^2\,c_{22}^k,\\
\vspace{2mm} \bar c_{12}^k = u_{\bar u}\,u_{\bar v} \,c_{11}^k  +
(u_{\bar u} \, v_{\bar v}+ u_{\bar v}\, v_{\bar u})\, c_{12}^k +
v_{\bar u} \, v_{\bar v} \, c_{22}^k,\\
\vspace{2mm} \bar c_{22}^k = u_{\bar v}^2\,c_{11}^k  + 2 u_{\bar
v} \, v_{\bar v}\, c_{12}^k + v_{\bar v}^2\,c_{22}^k.
\end{array} \quad  \quad  (k = 1,2)\leqno{(2.3)}$$

Using (2.3), we obtain
$$\begin{array}{l}
\vspace{2mm} \overline{\Delta}_1 = J \left(u_{\bar u}^2\,\Delta_1
+
u_{\bar u} \, v_{\bar u}\, \Delta_2 + v_{\bar u}^2\,\Delta_3 \right);\\
\vspace{2mm} \overline{\Delta}_2 = J \left(2 u_{\bar u}\,u_{\bar
v}\,\Delta_1 + (u_{\bar u}\,v_{\bar v} + u_{\bar v} \, v_{\bar
u})\,\Delta_2
+ 2 v_{\bar u}\, v_{\bar v}\,\Delta_3 \right);\\
\vspace{2mm} \overline{\Delta}_3 = J \left(u_{\bar v}^2\,\Delta_1
+ u_{\bar v} \, v_{\bar v}\, \Delta_2 + v_{\bar v}^2\,\Delta_3
\right).
\end{array}\leqno{(2.4)}$$

If $\bar E = g(z_{\bar u}, z_{\bar u})$, $\bar F = g(z_{\bar u},
z_{\bar v})$ and $\bar G = g(z_{\bar v}, z_{\bar v})$, then we
have
$$\begin{array}{l}
\vspace{2mm}
\bar E=u_{\bar u}^2\,E+2\,u_{\bar u}v_{\bar u}\,F+v_{\bar u}^2\,G,\\
\vspace{2mm} \bar F=u_{\bar u}u_{\bar v}\,E+(u_{\bar u}v_{\bar
v}+v_{\bar u}u_{\bar v})\,F
+v_{\bar u}v_{\bar v}\,G,\\
\vspace{2mm} \bar G=u_{\bar v}^2\,E+2\,u_{\bar v}v_{\bar
v}\,F+v_{\bar v}^2\,G\end{array}\leqno{(2.5)}$$
and
$$\bar E \bar G - \bar F^2=J^2\,(EG-F^2)$$
or
$$\bar W = \varepsilon J\,W, \quad \varepsilon = {\rm sign} \, J. \leqno{(2.6)}$$
Taking into account (2.4) and (2.6), we find
$$\begin{array}{l}
\vspace{2mm}
\bar L=\varepsilon (u_{\bar u}^2\,L+2\,u_{\bar u}v_{\bar u}\,M+v_{\bar u}^2\,N),\\
\vspace{2mm} \bar M = \varepsilon (u_{\bar u}u_{\bar
v}\,L+(u_{\bar u}v_{\bar v}+v_{\bar u}u_{\bar v})\,M
+v_{\bar u}v_{\bar v}\,N),\\
\vspace{2mm} \bar N=\varepsilon (u_{\bar v}^2\,L+2\,u_{\bar
v}v_{\bar v}\,M+v_{\bar v}^2\,N).
\end{array}\leqno{(2.7)}$$

Further we denote
$$\begin{array}{l}
\vspace{2mm}
\displaystyle{\gamma_1^1=\frac{FM-GL}{EG-F^2}, \qquad \gamma_1^2=\frac{FL-EM}{EG-F^2}},\\
\vspace{2mm} \displaystyle{\gamma_2^1=\frac{FN-GM}{EG-F^2}, \qquad
\gamma_2^2=\frac{FM-EN}{EG-F^2}}
\end{array}\leqno{(2.8)}$$
and consider the linear map
$$\gamma: T_pM^2 \rightarrow T_pM^2$$
determined by the conditions
$$\begin{array}{l}
\vspace{2mm}
\gamma(z_u)=\gamma_1^1z_u+\gamma_1^2z_v,\\
\vspace{2mm} \gamma(z_v)=\gamma_2^1z_u+\gamma_2^2z_v,
\end{array}
\qquad \gamma = \left(
\begin{array}{cc}
\gamma_1^1 & \gamma_1^2\\
[2mm] \gamma_2^1 & \gamma_2^2
\end{array} \right).
\leqno{(2.9)}$$
Then a tangent vector $X=\lambda z_u+ \mu z_v$ is transformed
into the vector $X'=\gamma(X)=\lambda' z_u+\mu' z_v$ so that
$$\left(
\begin{array}{l}
\lambda'\\
[2mm] \mu' \end{array}\right)=\gamma^t \left(
\begin{array}{l}
\lambda\\
[2mm] \mu \end{array} \right).$$

We have

\begin{lem}\label{L:2}
The linear map $\gamma$ given by $(2.9)$ is geometrically
determined.
\end{lem}

\emph{Proof:} Let the change of the parameters be given by (2.2).
Then we have
$$\left(
\begin{array}{l}
z_{\bar u}\\
[2mm] z_{\bar v}
\end{array} \right) = T \left(
\begin{array}{l}
z_u\\
[2mm] z_v
\end{array} \right), \qquad
T=\left(
\begin{array}{cc}
u_{\bar u} & v_{\bar u}\\
[2mm] u_{\bar v} & v_{\bar v}
\end{array} \right).$$
If we denote
$$g=\left(
\begin{array}{cc}
E & F\\
[2mm] F & G \end{array} \right), \qquad h=\left(
\begin{array}{cc}
L & M\\
[2mm] M & N \end{array} \right),$$ then the defining conditions
(2.8) imply $\gamma = -hg^{-1}$.

With respect to the new coordinates $(\bar u, \bar v)$ the linear
map $\bar \gamma$ is determined by the equality $\bar \gamma=
-\bar h \bar g^{-1}.$

On the other hand, the equalities (2.5) and (2.7) express that
$$\bar g=TgT^t, \qquad \bar h=\varepsilon \, ThT^t.$$
Thus we obtain $\bar \gamma = -\bar h \bar g^{-1}=\varepsilon \,
T\gamma T^{-1},$ which implies that $\bar \gamma =\varepsilon \,
\gamma$.

Further, let $\{\widetilde{e}_1, \widetilde{e}_2\}$ be another
orthonormal normal frame field of $M^2$.  Then
$$\begin{array}{l}
\vspace{2mm}
e_1 = \cos \theta \, \widetilde{e}_1 + \varepsilon' \sin \theta \, \widetilde{e}_2;\\
\vspace{2mm} e_2 = - \sin \theta \, \widetilde{e}_1 + \varepsilon'
\cos \theta \, \widetilde{e}_2;
\end{array}
\qquad \theta = \angle (\widetilde{e}_1, e_1),$$ and $\varepsilon'
= 1 \; (\varepsilon' =-1)$ if the normal frame fields $\{e_1, e_2\}$ and
$\{\widetilde e_1, \widetilde e_2\}$ have the same (opposite)
orientation. The relation between the corresponding functions
$c_{ij}^k$ and $\widetilde{c}_{ij}^k$, $i,j,k = 1,2$ is given by
the equalities
$$\begin{array}{l}
\vspace{2mm}
\widetilde{c}_{ij}^1 = \cos \theta \, c_{ij}^1 - \sin \theta \, c_{ij}^2;\\
\vspace{2mm} \widetilde{c}_{ij}^2 = \varepsilon' (\sin \theta \,
c_{ij}^1 + \cos \theta \, c_{ij}^2);
\end{array} \quad i,j = 1,2.$$

Thus,  $\widetilde{\Delta}_i = \varepsilon' \Delta_i$, $i =
1,2,3$, and $\widetilde{L} = \varepsilon' L, \; \widetilde{M} =
\varepsilon' M, \; \widetilde{N} = \varepsilon' N$, which imply that
$\widetilde \gamma= \varepsilon' \, \gamma $.

\qed

\vskip 2mm
The linear map $\gamma: T_pM^2 \rightarrow T_pM^2$ is said to
be the \emph{Weingarten map} at the point $p \in M^2$. The
following statement follows immediately from Lemma \ref{L:2}.

\begin{lem}\label{L:Lemma-invariant}
The functions
$$k := \det \gamma = \frac{LN - M^2}{EG - F^2}, \qquad
\varkappa :=-\frac{1}{2}\,{\rm tr}\, \gamma =\frac{EN+GL-2FM}{2(EG-F^2)}\leqno(2.10)$$
are invariants of the surface $M^2$.
\end{lem}

\vskip 2mm
It is clear that the sign of $\varkappa$ depends on the
orientations of the tangent plane and the normal space of $M^2$, while
$k$ is an absolute invariant.

The characteristic equation of the Weingarten map $\gamma$ in view
of Lemma \ref{L:Lemma-invariant} is
$$\nu^2 + 2 \varkappa \, \nu + k = 0.\leqno(2.11)$$
\vskip 2mm

If $X_1$ and $X_2$ are two tangent vectors at a point $p \in M^2$,
then $g(\gamma (X_1), X_2)=g(\gamma(X_2), X_1)$, i.e. $\gamma$ is
a symmetric linear operator and hence
$$\varkappa^2 - k \geq 0.\leqno(2.12)$$

Using the defining equalities (2.10), it follows that
$$4(\varkappa^2 - k)=\left(
\gamma_1^1-\gamma_2^2+2\frac{F}{E}\gamma_1^2\right)^2+4\frac{EG-F^2}{E^2}(\gamma_1^2)^2.$$
This equality implies that the condition $\varkappa^2 - k  = 0$ is equivalent
to the equalities $\gamma_1^1=\gamma_2^2, \; \gamma_1^2=0$, i.e.
to the conditions
$$L = \rho E, \quad M = \rho F, \quad N = \rho G, \qquad \rho\in {\R}.\leqno(2.13)$$

Thus we get the following equivalence at a point $p \in M^2$:
$$L = M = N = 0 \quad \iff \quad k = \varkappa = 0.\leqno(2.14)$$

As in the classical case (for a surface $M^2$ in ${\E}^3$), the invariants $k$ and
$\varkappa$ divide the points of $M^2$ into four types. A point $p \in M^2$
is said to be:
\vskip 2mm
\emph{flat}, \; if
\; $k = \varkappa = 0$; \vskip 1mm \emph{elliptic}, \; if \; $k > 0$; \vskip 1mm
\emph{parabolic}, \; if \; $k = 0$, \; $\varkappa \neq 0$; \vskip 1mm
\emph{hyperbolic}, \; if \; $k < 0$.

\vskip 2mm
Let $X=\lambda z_u+\mu z_v, \,\, (\lambda,\mu)\neq(0,0)$ be a tangent
vector at a point $p \in M^2$. The Weingarten map $\gamma$ determines a second
fundamental form of the surface $M^2$ at $p \in M^2$ as follows:
$$II(\lambda,\mu)=-g(\gamma(X),X)=L\lambda^2+2M\lambda\mu+N\mu^2, \quad \lambda, \mu \in {\R}.$$

\vskip 2mm First we study the class of surfaces whose points are
flat.

\section{Surfaces consisting of flat points} \label{S:flat}

In this section we consider surfaces $M^2: z = z(u,v), \, \, (u,v)
\in {\mathcal D}$ consisting of flat points, i.e. surfaces satisfying
the conditions
$$k(u,v)=0, \quad \varkappa(u,v)=0, \qquad (u,v) \in \mathcal D.\leqno(3.1)$$

We give a local geometric description of these surfaces.

For the sake of simplicity, we shall assume that the
parametrization of $M^2$ is orthogonal, i.e. $F=0$. Denote the unit vector
fields $\displaystyle{x=\frac{z_u}{\sqrt E}, \; y=\frac{z_v}{\sqrt
G}}$. Then we write (2.1) in the form
$$\begin{array}{l}
\vspace{2mm} \nabla'_xx=\quad \quad \quad
\gamma_1\,y+\;\;\displaystyle{\frac{c_{11}^1}{E}\;\;e_1
\;+\;\frac{c_{11}^2}{E}\;\;e_2},\\
\vspace{2mm} \nabla'_xy=-\gamma_1\,x\quad \quad \;
+\displaystyle{\frac{c_{12}^1}{\sqrt{EG}}\,e_1
+\frac{c_{12}^2}{\sqrt{EG}}\,e_2,}\\
\vspace{2mm} \nabla'_yx=\quad\quad -\gamma_2\,y \;\;
+\displaystyle{\frac{c_{12}^1}{\sqrt{EG}}\,e_1
+\frac{c_{12}^2}{\sqrt{EG}}\,e_2,}\\
\vspace{2mm} \nabla'_yy=\;\;\gamma_2\,y
\quad\quad\quad+\;\;\displaystyle{\frac{c_{22}^1}{G}\;\;e_1
\;+\;\;\frac{c_{22}^2}{G}\;\;e_2}.\\
\end{array}\leqno(3.2)$$

Obviously, the surface $M^2$ lies in a 2-plane if and only if $M^2$ is
totally geodesic, i.e. $c_{ij}^k=0, \; i,j,k = 1, 2.$

Now, let at least one of the coefficients $c_{ij}^k$ not be zero.
Then
$${\rm rank}\left(
\begin{array}{ccc}
c_{11}^1 & c_{12}^1 & c_{22}^1\\
[2mm] c_{11}^2 & c_{12}^2 & c_{22}^2
\end{array}\right)=1 $$
and the vectors $\sigma (x,x), \; \sigma(x,y), \; \sigma(y,y)$ are
collinear. Let $\{b,l\}$ be a normal frame field of $M^2$,
consisting of orthonormal vector fields, such that $b$ is
collinear with $\sigma (x,x), \; \sigma(x,y)$, and $\sigma(y,y)$.
It is clear that the normal frame field $\{b,l\}$ is invariant.
Then the derivative formulas of $M^2$ can be written as follows:
$$\begin{array}{ll}
\vspace{2mm}
\nabla'_xx=\quad \quad \quad \gamma_1\,y+\;\nu_1\,b,& \qquad \quad
\nabla'_xb=-\nu_1\,x-\lambda\,y \quad \, \quad+\;\beta_1\,l,\\
\vspace{2mm}
\nabla'_xy=-\gamma_1\,x\quad \quad \; + \; \lambda\,b, & \qquad \quad
\nabla'_yb=-\lambda \, x -\nu_2\,y\quad \quad \,+ \beta_2 \,l,\\
\vspace{2mm}
\nabla'_yx=\quad\quad -\gamma_2\,y \; + \;\lambda\, b, & \qquad \quad
\nabla'_xl=\quad\quad \quad \quad \quad-\beta_1\,b,\\
\vspace{2mm}
\nabla'_yy=\;\;\gamma_2\,x \quad\quad \;\;+\;\nu_2\,b, & \qquad \quad
\nabla'_yl=\quad\quad \quad \quad \quad-\beta_2\,b,
\end{array}\leqno(3.3)$$
for some functions $\nu_1, \, \nu_2,\, \lambda,\,\beta_1,\, \beta_2, \gamma_1,\, \gamma_2$
on $M^2$.

The Gauss curvature $K$ of $M^2$ is expressed by
$$K = \nu_1\,\nu_2-\lambda^2.\leqno(3.4)$$

Further we denote $\beta = \beta_1^2+\beta_2^2$. It
follows immediately that $\beta$ does not depend on the change (2.2)
of the parameters.

Since the curvature tensor $R'$ of the connection $\nabla'$ is zero, then the equalities
$R'(x,y,b) = 0$ and $R'(x,y,l) = 0$ together with  (3.3) imply that
either $K = 0$ or $\beta = 0$.

\vskip 2mm
A surface $M^2$ is said to be \emph{planar} if there exists a
hyperplane $\E^3 \subset \E^4$ containing $M^2$. First we shall characterize the
planar surfaces.

\begin{prop}\label{P:3.1}
A surface $M^2$ is planar if and only if
$$k = 0, \quad \varkappa = 0, \quad \beta = 0.$$
\end{prop}

\emph{Proof:} I. Let $M^2 \subset {\E}^3$ and $b$ be the usual
normal to $M^2$ in ${\E}^3$. Choosing $l$ to be the normal to the
hyperplane ${\E}^3$, from (3.3) we get $L=M=N=0$ and $\beta=0$.

II.
Under the conditions $k = \varkappa = \beta = 0$, from (3.3) it follows that
$l={\rm const}$ and $M^2$ lies in a hyperplane ${\E}^3$ orthogonal
to $l$.
\qed

\vskip 3mm
A \textit{ruled surface} $M^2$ is a one-parameter system $\{g(v)\}, v \in J$
of straight lines $g(v)$, defined in an interval $J \subset \R$. The straight lines
$g(v)$ are called generators of $M^2$. A ruled surface $M^2 = \{g(v)\}, v \in J$ is
said to be {\it developable}, if the tangent space $T_pM^2$ at all
regular points $p$ of an arbitrary fixed generator $g(v)$ is one
and the same.

Each ruled surface $M^2$ can be parameterized as follows:
$$z(u,v) = x(v) + u \,e(v), \quad u \in \R, \,\, v \in J, \leqno(3.5)$$
where $x(v)$ and $e(v)$ are vector-valued functions, defined in
$J$, such that the vectors $e(v)$ and $x'(v) + u \,e'(v)$ are
linearly independent for all $v \in J$. The tangent space of $M^2$
is spanned by the vectors
$$\begin{array}{l}
\vspace{1mm}
z_u = e(v);\\
z_v = x'(v) + u\,e'(v).
\end{array}$$

The ruled surface $M^2$ determined by (3.5) is developable if and only
if the vectors $e(v)$, $e'(v)$ and $x'(v)$ are linearly dependent.

We shall characterize the developable ruled surfaces in terms of
the invariants $k$, $\varkappa$ and the Gauss curvature $K$.

\begin{prop}\label{P:3.2}
A surface $M^2$ is locally a developable ruled surface if and only if
$$k = 0, \quad \varkappa = 0, \quad K = 0.$$
\end{prop}

\emph{Proof:} I. Let $M^2$ be a developable ruled surface, defined
by the equality (3.5), where $e(v)$, $e'(v)$ and $x'(v)$ are linearly
dependent. Without loss of generality we assume that $e^2(v) = 1$.
Then, the vector fields $e(v)$ and $e'(v)$ are orthogonal and the
tangent space of $M^2$ is $\span \{e(v), e'(v)\}$. Since $x'(v)
\in \span \{e(v), e'(v)\}$, then $x'(v)$ is decomposed in the form
$x'(v) = p(v)\, e(v) + q(v)\, e'(v)$ for some functions $p(v)$ and
$q(v)$. Hence, the tangent space of $M^2$ is spanned by
$$\begin{array}{l}
\vspace{1mm}
z_u = e;\\
z_v = p \,e + (u + q)\,e'.
\end{array}$$

Considering only the regular points of $M^2$ (where $u \neq -q$),
we choose an orthonormal tangent frame field $\{x, y\}$ of $M^2$
in the following way:
$$\begin{array}{l}
\vspace{1mm}
x = e = z_u;\\
y = \displaystyle{\frac{e'}{\sqrt{(e')^2}} = - \frac{p}{(u +
q)\sqrt{(e')^2}}\, z_u + \frac{1}{(u + q)\sqrt{(e')^2}} \, z_v}.
\end{array} \leqno{(3.6)}$$

Since the tangent space of $M^2$ does not depend on the parameter
$u$, then the normal space of $M^2$ is spanned by vector fields
$b_1(v),b_2(v)$. With respect to the basis $\{e(v), e'(v),$
$b_1(v), b_2(v)\}$ the derivatives of $b_1(v)$ and $b_2(v)$ are
decomposed in the form
$$\begin{array}{l}
\vspace{1mm}
b_1' = - c_1\,e' + c_0 \,b_2,\\
b_2' = - c_2\,e' - c_0 \,b_1,\\
\end{array}\leqno{(3.7)}$$ where
$c_0$, $c_1$, $c_2$ are functions of $v$.

Then the equalities (3.6) and (3.7) imply
$$\begin{array}{l}
\vspace{1mm}
\nabla'_x b_1 = 0,\\
\vspace{1mm} \nabla'_y b_1 = \displaystyle{-\frac{c_1}{u +q}\,
y + \frac{c_0}{(u +q)\sqrt{(e')^2}} \, b_2},\\
\vspace{1mm}
\nabla'_x b_2 = 0,\\
\vspace{1mm} \nabla'_y b_2 = \displaystyle{-\frac{c_2}{u +q}\, y -
\frac{c_0}{(u +q)\sqrt{(e')^2}} \, b_1}.
\end{array}$$
Consequently, $L = M = N = 0$ and $K = 0$.

II. Let $M^2$ be a surface for which $L = M = N = 0$ and $K = 0$.
We consider an orthonormal frame field $\{x, y, b, l\}$ of $M^2$,
satisfying the equalities (3.3). Since $K = 0$, then
$\nu_1\,\nu_2-\lambda^2 = 0$. If $\nu_1 = \nu_2 = 0$, then $M^2$
lies in a plane $\E^2$. So we assume that there exists a
neighborhood $\widetilde{\mathcal D} \subset \mathcal D$ such that
${\nu_2}_{|\widetilde{\mathcal D}} \neq 0$ (or
${\nu_1}_{|\widetilde{\mathcal D}} \neq 0$) and we consider the
surface $\widetilde{M}^2 = M^2_{|\widetilde{\mathcal D}}$.

Let $\{\overline{x},\overline{y}\}$ be the orthonormal tangent
frame field of $\widetilde{M}^2$, defined by
$$\begin{array}{l}
\vspace{2mm}
\overline{x} = \cos \varphi \, x + \sin \varphi \, y;\\
\vspace{2mm} \overline{y} = - \sin \varphi \, x +  \cos \varphi \,y,
\end{array}$$
where $\tan \varphi = \displaystyle{- \frac{\lambda}{\nu_2}}$.
Then $\sigma(\overline{x}, \overline{x}) = 0$, \,
$\sigma(\overline{x}, \overline{y}) = 0$. So the formulas (3.3) take the form
$$\begin{array}{ll}
\vspace{2mm}
\nabla'_{\overline{x}}\,\overline{x} = \quad \quad \quad \overline{\gamma}_1\,\overline{y},&
\qquad \quad \nabla'_{\overline{x}}\,b = \quad \quad \quad \quad \quad \quad \,\,\overline{\beta}_1\,l,\\
\vspace{2mm}
\nabla'_{\overline{x}}\,\overline{y} = -\overline{\gamma}_1\,\overline{x},&\qquad \quad
\nabla'_{\overline{y}}\,b = - \overline{\nu}_2\, \overline{y} \quad\quad \quad +\overline{\beta}_2\,l,\\
\vspace{2mm}
\nabla'_{\overline{y}}\,\overline{x} = \quad \quad \quad - \overline{\gamma}_2\,\overline{y},&
\qquad \quad\nabla'_{\overline{x}}\,l = \quad\quad \quad -\overline{\beta}_1\,b,\\
\vspace{2mm}
\nabla'_{\overline{y}}\,\overline{y} = \overline{\gamma}_2\,\overline{x} \quad \quad \quad \quad
+ \overline{\nu}_2\, b,& \qquad \quad
\nabla'_{\overline{y}}\,l = \quad\quad \quad -\overline{\beta}_2\,b,
\end{array}$$
where $\overline{\nu}_2 \neq 0$.

Since the curvature tensor $R'$ is zero, then the equalities
$R'(\overline{x},\overline{y},b) = 0$ and $R'(\overline{x},\overline{y},l) = 0$ imply that
$$\overline{\gamma}_1 = 0, \quad \quad \overline{\beta}_1 = 0.$$
Hence,
$$\begin{array}{ll}
\vspace{2mm}
\nabla'_{\overline{x}}\,\overline{x} = 0, & \quad \quad \nabla'_{\overline{x}}\,b = 0,\\
\vspace{2mm} \nabla'_{\overline{x}}\,\overline{y} = 0, & \quad
\quad \nabla'_{\overline{x}}\,l = 0.
\end{array}$$

Let $p= z(\overline{u}_0,\overline{v}_0), \,\,
(\overline{u}_0,\overline{v}_0) \in \widetilde{\mathcal D}$ be an
arbitrary point of  $\widetilde{M}^2$ and $c_1: z(\overline{u}) =
z(\overline{u},\overline{v}_0)$ be the integral curve of the
vector field $\overline{x}$, passing through $p$. From
$\nabla'_{\overline{x}}\,\overline{x} = 0$ it follows that $c_1$
is a straight line. Hence, for each point $p \in \widetilde{M}^2$
there exists a straight line passing through $p$, i.e.
$\widetilde{M}^2$ is a ruled surface. Moreover, since
$\nabla'_{\overline{x}}\,b = 0$ and $\nabla'_{\overline{x}}\,l =
0$ then the normal space $\span \{b, l\}$ of
$\widetilde{M}^2$ is constant at the points of $c_1$ and
hence, the tangent space $\span\{\overline{x},\overline{y}\}$ of
$\widetilde{M}^2$ at the points of $c_1$ is one and the same.
Consequently, $\widetilde{M}^2$ is a developable surface.
\qed

\vskip 2mm
From now on we exclude the flat points from our considerations.

\section{Minimal surfaces} \label{S:minimal}

We recall that a surface $M^2$ is said to be \textit{minimal} if the mean curvature vector
$\displaystyle{\frac{1}{2}\, \tr \, \sigma =0}$.
In this section we characterize the minimal surfaces in terms of the invariants
$k$ and $\varkappa$.

\vskip 2mm
\begin{prop}\label{P:minimal}
A surface $M^2$ in $\E^4$ is minimal if and only if
$$\varkappa^2 - k = 0.$$
\end{prop}

\vskip 2mm
\emph{Proof:} Without loss of generality we assume that $F = 0$ and denote the unit vector fields
$\displaystyle{x=\frac{z_u}{\sqrt E}, \; y=\frac{z_v}{\sqrt G}}$.
Then we have
$$\begin{array}{l}
\vspace{2mm} \nabla'_xx=\quad \quad \quad
\gamma_1\,y+\;\;\displaystyle{\frac{c_{11}^1}{E}\;\;e_1
\;+\;\frac{c_{11}^2}{E}\;\;e_2},\\
\vspace{2mm} \nabla'_xy=-\gamma_1\,x\quad \quad \;
+\displaystyle{\frac{c_{12}^1}{\sqrt{EG}}\,e_1
+\frac{c_{12}^2}{\sqrt{EG}}\,e_2,}\\
\vspace{2mm} \nabla'_yx=\quad\quad -\gamma_2\,y \;\;
+\displaystyle{\frac{c_{12}^1}{\sqrt{EG}}\,e_1
+\frac{c_{12}^2}{\sqrt{EG}}\,e_2,}\\
\vspace{2mm} \nabla'_yy=\;\;\gamma_2\,y
\quad\quad\quad+\;\;\displaystyle{\frac{c_{22}^1}{G}\;\;e_1
\;+\;\;\frac{c_{22}^2}{G}\;\;e_2}.\\
\end{array}$$

I. Let $\tr \, \sigma = \sigma(x, x)+\sigma(y,y)=0$. Then
$$\Delta_2 = \left|%
\begin{array}{cc}
\vspace{2mm}
  c_{11}^1 & c_{22}^1 \\
  c_{11}^2 & c_{22}^2 \\
\end{array}%
\right|=0, \quad \quad \frac{\Delta_3}{G}=\frac{\Delta_1}{E},$$
and hence
$$L = \rho E, \quad M = \rho F, \quad N = \rho G.$$
The last equalities imply that $\varkappa^2 - k = 0$.

II. Let $\varkappa^2 - k = 0$. Then
$$L = \rho E, \quad M = \rho F, \quad N = \rho G.$$
The condition $F=0$ implies that $M=0$. Then
$\displaystyle{\left|%
\begin{array}{cc}
\vspace{2mm}
  c_{11}^1 & c_{22}^1 \\
  c_{11}^2 & c_{22}^2 \\
\end{array}%
\right|=0}$
and $c_{22}^1=\widetilde{\rho} c_{11}^1, \; c_{22}^2=\widetilde{\rho}
c_{11}^2$.
Further, the equality
$\displaystyle{\frac{L}{E}=\frac{N}{G}}$ implies that
$\displaystyle{\widetilde{\rho}=-\frac{G}{E}}$. Hence $\tr \, \sigma = 0$.
\qed
\vskip 2mm

\section{Surfaces of general type} \label{S:general}
From now on we consider surfaces, satisfying the condition
$$\varkappa^2 - k \neq 0$$
and call them surfaces  \textit{of general type}.

As in the classical differential geometry of surfaces in
$\E^3$ the second fundamental form  determines conjugate tangents at a point $p$ of $M^2$.
A tangent $g: X=\lambda z_u+\mu z_v$ is said to be
\emph{principal} if it is perpendicular to its conjugate. The
equation for the principal tangents at a point $p \in M^2$ is
$$\left|\begin{array}{cc}
E & F\\
[2mm] L & M \end{array}\right| \lambda^2+ \left|\begin{array}{cc}
E & G\\
[2mm] L & N \end{array}\right| \lambda \mu+
\left|\begin{array}{cc}
F & G\\
[2mm] M & N \end{array}\right| \mu^2=0.$$

A line $c: u=u(q), \; v=v(q); \; q\in J$ on $M^2$ is said to be a
\emph{principal} curve if its tangent at any point is principal.

The surface $M^2$ is parameterized with respect to the
principle lines if and only if
$$F=0, \qquad M=0.$$

Let $M^2$ be parameterized with respect to the principal lines and
denote the unit vector fields
$\displaystyle{x=\frac{z_u}{\sqrt E}, \; y=\frac{z_v}{\sqrt G}}$.

Since the mean curvature vector field
$H \neq 0$, we determine
the unit normal vector field $b$ by the equality
$b=\displaystyle{\frac{H}{\Vert H \Vert}}$.
Further we denote by $l$ the unit normal vector field such that
$\{x,y,b,l\}$ is a positive oriented orthonormal frame field of
$M^2$.
Thus we obtain a geometrically determined orthonormal frame field $\{x,y,b,l\}$
at each point $p \in M^2$.
With respect to the frame field $\{x,y,b,l\}$ we have the
following Frenet-type derivative formulas:
$$\begin{array}{ll}
\vspace{2mm} \nabla'_xx=\quad \quad \quad \gamma_1\,y+\,\nu_1\,b;
& \qquad
\nabla'_xb=-\nu_1\,x-\lambda\,y\quad\quad \quad \,+\beta_1\,l;\\
\vspace{2mm} \nabla'_xy=-\gamma_1\,x\quad \quad \; + \; \lambda\,b
\; + \mu\,l;  & \qquad
\nabla'_yb=-\lambda\,x - \; \nu_2\,y\quad\quad \quad +\beta_2\,l;\\
\vspace{2mm} \nabla'_yx=\quad\quad \;-\gamma_2\,y \; + \lambda\,b
\; +\mu\,l;  & \qquad
\nabla'_xl= \quad \quad \quad \;-\mu\,y-\beta_1\,b;\\
\vspace{2mm} \nabla'_yy=\;\;\gamma_2\,x \quad\quad\quad+\nu_2\,b;
& \qquad \nabla'_yl=-\mu\,x \quad \quad \quad \;-\beta_2\,b,
\end{array}\leqno{(5.1)}$$
where  $\gamma_1 = - y(\ln \sqrt{E}), \,\, \gamma_2 = - x(\ln
\sqrt{G})$ and $\mu \neq 0$.

Hence we have
$$k = - 4\nu_1\,\nu_2\,\mu^2, \quad \quad \varkappa = (\nu_1-\nu_2)\mu, \quad \quad
K = \nu_1\,\nu_2- (\lambda^2 + \mu^2).\leqno(5.2)$$

\vskip 3mm
\noindent
\textbf{Remark 1.}
We note that we determine the tangent frame field $\{x,y\}$ by the Weingarten map
(the second fundamental form $II$) and the normal frame field $\{b,l\}$ -
by the mean curvature vector field, while the Frenet-cross section in the sense of
\={O}tsuki  diagonalizes a quadratic form in the normal space. In general
the geometric frame field $\{b,l\}$ is not a Frenet-cross section. Finding the relation between
$\{b,l\}$  and the Frenet-cross section of \={O}tsuki we derive the following relation between
the invariant $k$ and the curvatures $\lambda_1$ and $\lambda_2$ of \={O}tsuki:
$$k = 4 \lambda_1 \lambda_2.$$
The same formula is valid in the cases of minimal surfaces and surfaces consisting of
flat points.

\vskip 3mm
Using (5.1) we find the length of the mean curvature vector $\Vert H \Vert$ and
taking into account (5.2) we obtain the formula
$$\Vert H \Vert = \displaystyle{\frac{\sqrt{\varkappa^2-k}}{2 |\mu |}},$$
which shows that $|\mu|$ is expressed by the invariants $k$, $\varkappa$ and
the mean curvature function.

Let $z = g(z,x)\,x + g(z,y)\,y$ be an arbitrary tangent vector field of $M^2$.
We define the one-form $\theta$ by the equality
$$\theta (z) = g(\nabla'_z b,l).$$
Then the formulas (5.1) imply that
$$\theta (z)= g(\beta_1 \,x + \beta_2 \,y,z),$$
which shows that the one-form $\theta$ corresponds to the tangent vector field
$\beta_1 \,x + \beta_2 \,y$ and
$$\Vert \theta \Vert = \sqrt{\beta_1^2 + \beta_2^2}.$$

Using that $R'(x,y,x) = 0$, $R'(x,y,y) = 0$, $R'(x,y,b) =
0$ and $R'(x,y,l) = 0$, we get the following integrability
conditions:
$$\begin{array}{l}
\vspace{2mm}
\nu_1 \,\nu_2 - (\lambda^2 + \mu^2) = x(\gamma_2) + y(\gamma_1) - \left((\gamma_1)^2 + (\gamma_2)^2\right);\\
\vspace{2mm}
2\mu\, \gamma_2 + \nu_1\,\beta_2 - \lambda\,\beta_1 = x(\mu);\\
\vspace{2mm}
2\mu\, \gamma_1 - \lambda\,\beta_2 + \nu_2\,\beta_1 = y(\mu);\\
\vspace{2mm}
2\lambda\, \gamma_2 + \mu\,\beta_1 - (\nu_1 - \nu_2)\,\gamma_1 = x(\lambda) - y(\nu_1);\\
\vspace{2mm}
2\lambda\, \gamma_1 + \mu\,\beta_2 + (\nu_1 - \nu_2)\,\gamma_2 = - x(\nu_2) + y(\lambda);\\
\gamma_1\,\beta_1 - \gamma_2\,\beta_2 + (\nu_1 - \nu_2)\,\mu  = -
x(\beta_2) + y(\beta_1).
\end{array}\leqno{(5.3)}$$

\vskip 3mm
At the end of this section we shall characterize the surfaces with flat normal connection
in terms of the invariant $\varkappa$.

\vskip 2mm
A surface $M^2$ is said to be of \textit{flat normal connection} \cite{Chen1}
if the normal curvature $R^{\bot}$ of $M^2$ is zero. The
equalities (5.1) imply that the normal curvature $R^{\bot}$ of
$M^2$ is expressed as follows:
$$\begin{array}{l}
\vspace{2mm}
R^{\bot}_b(x,y) = D_xD_yb - D_yD_xb - D_{[x,y]}b
= \left(x(\beta_2) - y(\beta_1) + \gamma_1\,\beta_1 - \gamma_2 \,\beta_2 \right)\,l,\\
\vspace{2mm} R^{\bot}_l(x,y) = D_xD_yl - D_yD_xl - D_{[x,y]}l = -
\left(x(\beta_2) - y(\beta_1) + \gamma_1\,\beta_1 - \gamma_2
\,\beta_2 \right)\,b.
\end{array}\leqno{(5.4)}$$

Taking in mind (5.4) and the last equality of (5.3) we get:
$$\begin{array}{l}
\vspace{2mm}
R^{\bot}_b(x,y) = - \varkappa \,l,\\
\vspace{2mm} R^{\bot}_l(x,y) = \varkappa \,b,
\end{array}$$
i.e.
$$\varkappa = g(R^{\bot}_l(x,y), b) = g(R^{\bot}(x,y)l, b).$$
Hence, the invariant $\varkappa$ is the curvature of the normal connection of $M^2$.

\vskip 2mm Thus the surfaces with flat normal connection are
characterized by the following

\begin{prop}\label{P:flat-normal-1}
A surface $M^2$ in $\E^4$ is of flat normal connection if and only if
$$\varkappa = 0.$$
\end{prop}

Obviously, $M^2$ is a surface with flat normal connection if and only if $\nu_1
= \nu_2 =:\nu$. So, the Frenet-type formulas (5.1) of a surface
$M^2$ with flat normal connection take the form:
$$\begin{array}{ll}
\vspace{2mm} \nabla'_xx=\quad \quad \quad \gamma_1\,y+\,\nu\,b;  &
\qquad
\nabla'_xb=-\nu\,x-\lambda\,y\quad\quad \quad \, +\beta_1\,l;\\
\vspace{2mm} \nabla'_xy=-\gamma_1\,x\quad \quad \; + \; \lambda\,b
\; + \mu\,l;  & \qquad
\nabla'_yb=-\lambda\,x - \; \nu\,y\quad\quad \quad +\beta_2\,l;\\
\vspace{2mm} \nabla'_yx=\quad\quad \;-\gamma_2\,y \; + \lambda\,b
\; +\mu\,l;  & \qquad
\nabla'_xl= \quad \quad \quad -\mu\,y-\beta_1\,b;\\
\vspace{2mm} \nabla'_yy=\;\;\gamma_2\,x \quad\quad\quad+\nu\,b; &
\qquad \nabla'_yl=-\mu\,x \quad \quad \quad -\beta_2\,b.
\end{array}\leqno{(5.5)}$$

Hence the invariants $k$ and $K$ are expressed by
$$k = - 4\nu^2\,\mu^2, \quad \quad K = \nu^2 - (\lambda^2 + \mu^2).$$

\vskip 3mm
\noindent
\textbf{Remark 2.}
The curvature of the normal connection of a surface $M^2$ in $\E^4$ is the Gauss
torsion $\varkappa_G$ of $M^2$, which is defined in terms of the of the ellipse of
normal curvature at a point $p \in M^2$ \cite{Aminov}.

\section{Rotational surfaces}\label{S:rotational}

Now we shall apply our theory to the class of the rotational surfaces in $\E^4$.

We denote by $Oe_1e_2e_3$ a fixed orthonormal base of ${\E}^3$.
Let $\textbf{\textit{c}}: \widetilde{z} = \widetilde{z}(u),
\,\, u \in J$ be a smooth curve in $\E^3$, parameterized by
$$\widetilde{z}(u) = \left( x_1(u), x_2(u), r(u)\right); \quad u \in J.$$
We denote by $\textbf{\textit{c}}_1$ the projection of $\textbf{\textit{c}}$ on
the 2-dimensional plane $Oe_1e_2$.

Without loss of generality we can assume that
$\textbf{\textit{c}}$ is parameterized with respect to the
arc-length, i.e. $(x_1')^2 + (x_2')^2 + (r')^2 = 1$. We assume also
that $r(u)>0,\,\, u \in J$. Let us consider the rotational surface
$M^2$ in $\E^4$ given by
$$z(u,v) = \left( x_1(u), x_2(u), r(u) \cos v, r(u) \sin v\right);
\quad u \in J,\,\,  v \in [0; 2\pi).\leqno{(6.1)}$$

The tangent space of $M^2$ is spanned by the vector fields
$$\begin{array}{l}
\vspace{2mm}
z_u = \left(x_1', x_2', r' \cos v, r' \sin v \right);\\
\vspace{2mm} z_v = \left(0, 0, -r \sin v, r \cos v \right).
\end{array}$$
Hence,
$$E = 1; \quad F = 0; \quad G = r^2(u); \quad W = r(u).$$
We consider the following orthonormal tangent vector fields
$$\begin{array}{l}
\vspace{2mm}
\overline{x} = \left(x_1', x_2', r' \cos v, r' \sin v \right);\\
\vspace{2mm}
\overline{y} = \left(0, 0, - \sin v, \cos v \right),
\end{array}$$
i.e. $z_u = \overline{x}; \,\, z_v = r\, \overline{y}$. The second
partial derivatives of $z(u,v)$ are expressed as follows
$$\begin{array}{l}
\vspace{2mm}
z_{uu} = \left(x_1'', x_2'', r'' \cos v, r'' \sin v \right);\\
\vspace{2mm}
z_{uv} = \left(0, 0, -r' \sin v, r' \cos v \right);\\
\vspace{2mm} z_{vv} = \left(0, 0, -r \cos v, -r \sin v \right).
\end{array}$$

Let $\kappa$ and $\tau$ be the curvature and the torsion of the
curve $\textbf{\textit{c}}$ (considered as a curve in $\E^3$). We
consider the normal vector fields $e_1$ and $e_2$, defined by
$$\begin{array}{l}
\vspace{2mm}
e_1 = \displaystyle{\frac{1}{\kappa}\left(x_1'', x_2'', r'' \cos v, r'' \sin v \right)};\\
\vspace{2mm} e_2 = \displaystyle{\frac{1}{\kappa} \left( x_2' r''
- x_2'' r' , x_1'' r' - x_1' r'', (x_1' x_2'' - x_1'' x_2') \cos v, (x_1' x_2'' - x_1''
x_2') \sin v \right)}.
\end{array}$$
Now it is easy to calculate that
$$L = 0; \quad M = - (x_1' x_2'' - x_1'' x_2'); \quad N = 0.$$
Hence,
$$k = - \displaystyle{\frac{(x_1' x_2'' - x_1'' x_2')^2}{r^2}}; \quad \quad \varkappa = 0.$$

Applying Proposition \ref{P:flat-normal-1} we get
\begin{cor}\label{C:rotational}
Any rotational surface $M^2$ in $\E^4$, defined by $(6.1)$, is a
surface with flat normal connection.
\end{cor}

Let us denote the curvature of the plane curve $\textbf{\textit{c}}_1$ by
$\kappa_1 = x_1' x_2'' - x_1'' x_2'$. Then with respect
to the frame field $\{\overline{x}, \overline{y}, e_1, e_2\}$ the
derivative formulas of $M^2$ look like:
$$\begin{array}{ll}
\vspace{2mm} \nabla'_{\overline{x}}\overline{x} = \quad \quad
\quad \quad \quad  \kappa\,e_1 ;  & \qquad
\nabla'_{\overline{x}}e_1 = - \kappa\,\overline{x}\quad
\quad \quad \quad \quad \quad +\tau\,e_2;\\
\vspace{2mm} \nabla'_{\overline{x}}\overline{y} = 0;  & \qquad
\nabla'_{\overline{y}} e_1 = \quad  \quad \quad
\displaystyle{\frac{r''}{\kappa\,r}}\,\overline{y} ;\\
\vspace{2mm} \nabla'_{\overline{y}}\overline{x} = \quad
\quad\;\displaystyle{\frac{r'}{r}}\,\overline{y} ;  & \qquad
\nabla'_{\overline{x}}e_2 = \quad \quad \quad \quad \quad\quad -\tau\,e_1;\\
\vspace{2mm} \nabla'_{\overline{y}}\overline{y} =
-\displaystyle{\frac{r'}{r}}\,\overline{x} \quad \;-
\displaystyle{\frac{r''}{\kappa\,r}}\,e_1 -
\displaystyle{\frac{\kappa_1}{\kappa\,r}}\,e_2; & \qquad
\nabla'_{\overline{y}}e_2 = \quad \quad  \quad
\displaystyle{\frac{\kappa_1}{\kappa\,r}}\,\overline{y}.
\end{array}$$
So, the Gauss curvature of $M^2$ is:
$$K = - \displaystyle{\frac{r''}{r}}.$$

Obviously $M^2$ is not parameterized with respect to the principal lines.
The principal tangents of $M^2$ are:
$$\begin{array}{l}
\vspace{2mm}
x = \displaystyle{\frac{\sqrt{2}}{2}\, \overline{x} + \frac{\sqrt{2}}{2} \, \overline{y}};\\
\vspace{2mm} y = \displaystyle{\frac{\sqrt{2}}{2}\, \overline{x} -
\frac{\sqrt{2}}{2} \, \overline{y}}.
\end{array}$$

With respect to the geometric frame field $\{x,y,b,l\}$ the Frenet-type
formulas (5.5) hold good, where
$$\begin{array}{ll}
\vspace{2mm} \gamma_1 = \gamma_2 =
-\displaystyle{\frac{\sqrt{2}}{2}\,\frac{r'}{r}};& \quad \quad
\nu = \displaystyle{\frac{\sqrt{(\kappa^2 \,r - r'')^2 + (\kappa_1)^2}}{2\kappa \,r}};\\
\vspace{2mm} \lambda =\displaystyle{\frac{\kappa^4 \,r^2 -
(r'')^2 - (\kappa_1)^2}{2\kappa\,r}\sqrt{(\kappa^2 \,r -
r'')^2 + (\kappa_1)^2}}; & \quad\quad \mu
=\displaystyle{\frac{\kappa\,\kappa_1}{\sqrt{(\kappa^2
\,r - r'')^2 + (\kappa_1)^2}}}.
\end{array}$$
Consequently, the invariants $k$, $\varkappa$ and $K$ of the rotational surface $M^2$
are:
$$k = - \displaystyle{\frac{(\varkappa_1)^2}{r^2}}; \quad \quad \varkappa = 0; \quad \quad
K = - \displaystyle{\frac{r''}{r}}.$$

\vskip 3mm
At the end of the section we shall describe all rotational surfaces, for which the invariant
$k$ is constant.

\vskip 1mm
1. The invariant $k = 0$ if and only if $\kappa_1 = 0$,
which means that the projection of the curve $\textbf{\textit{c}}$ on the plane $Oe_1e_2$
is a straight line. There are two subcases:

\hskip 10mm
1.1. If $K = 0$, i.e. $r'' = 0$,  then $M^2$ is a developable ruled surface.

\hskip 10mm
1.2. If  $K \neq 0$, i.e. $r'' \neq 0$, then $M^2$ is a planar surface.

\vskip 1mm
2. The invariant $k = {\rm const}$ ($k \neq 0$) if and only if
$r(u) = a\, (x_1' x_2'' - x_1'' x_2')$,  $a = {\rm const}$.
Moreover, if $r(u)$ satisfies $r''(u) = c\, r(u)$, then the Gauss curvature $K$
is also a constant.

\end{document}